\documentclass{ifacconf}
\usepackage{graphicx}      % include this line if your document contains figures
\usepackage{natbib}        % required for bibliography
\usepackage{amsmath}
\usepackage{amssymb}
\usepackage{amsfonts}
\usepackage{enumitem}
\usepackage{comment}

\DeclareMathOperator{\sign}{sgn}
\newcommand{\R}{\mathbb{R}}
\renewcommand{\S}{\mathbb{S}}
\newcommand{\N}{\mathbb{N}}
\newcommand{\bmat}[1]{\begin{bmatrix}#1\end{bmatrix}}
\newcommand{\norm}[1]{\lVert{#1}\rVert}

\newcommand{\mcl}[1]{\mathcal{#1}}

\begin{document}
\begin{frontmatter}
\title{Bounding the Settling Time of Finite-Time Stable Systems using Sum of Squares}

\author{Sengiyumva Kisole, Kunal Garg, Matthew Peet}

\address{The School for the Engineering of Matter, Transport and Energy, Arizona State University, Tempe, AZ 85298, USA \\
(Emails: sengi.kisole@asu.edu, kgarg24@asu.edu, mpeet@asu.edu)}

\begin{abstract}                % Abstract of 50--100 words
Finite-time stability (FTS) of a differential equation guarantees that solutions reach a given equilibrium point in finite time, where the time of convergence depends on the initial state of the system. For traditional stability notions such as exponential stability, the convex optimization framework of Sum-of-Squares (SoS) enables computation of polynomial Lyapunov functions to certify stability. However, finite-time stable systems are characterized by non-Lipschitz, non-polynomial vector fields, rendering standard SoS methods inapplicable. To this end, we show that computation of a non-polynomial Lyapunov function certifying finite-time stability can be reformulated as feasibility of a set of polynomial inequalities under a particular transformation. As a result, SoS can be utilized to verify FTS and obtain a bound on the settling time. Numerical examples are used to demonstrate the accuracy of the conditions in both certifying finite-time stability and bounding the settling time. %This work represents the first combination of SoS programming with settling time bounds for finite-time stable systems.
\end{abstract}\vspace{-2mm}
\begin{keyword}
Finite-Time Stability, Sum of Squares Programming, Lyapunov Functions
\end{keyword}                    
\end{frontmatter}
%====================================================================================================================================================================
\section{Introduction}\vspace{-2mm}
For nonlinear systems, there exist multiple non-equivalent notions of stability, including: Lyapunov stability, asymptotic stability, exponential stability, and rational stability %(See~\citep{khalil1991control} and~\citep{bacciotti2005liapunov}). 
However, all of these notions of stability are asymptotic in the sense that none of them imply that solutions of a given system will ever reach a given stable set. By contrast, the notion of \textit{finite-time stability} of a set, $X$, ensures that for any initial condition, $x$, there exists a time, $T(x)$ such that if $x(t)$ is a solution then $x(t)\in X$ for all $t \ge T(x)$.

Finite-time stability plays a critical role in control applications, such as the design of robust controllers for sliding mode control (SMC) --- see \citep{polyakov2012unified, polyakov2014stability}. %Specifically, in SMC, asymptotic stability is only guaranteed once the solution reaches the sliding surface and hence twisting \citep{torres2017design, levant1993sliding} and super-twisting algorithms \citep{basin2014supertwisting, seeber2018lyapunov} are explicitly designed to achieve finite-time convergence to this surface. 
In addition, continuity of, and bounds on, the settling-time function $T(x)$ can be used, e.g. in robotics, to schedule a sequence of predefined motions -- each of which presumes some initial pose --- see~\citet{galicki2015finite} and~\citet{kong2020fuzzy}.
%. For example,~\citet{galicki2015finite} designs a class of robust non-singular terminal sliding mode controllers that guarantee finite time convergence of multi-DOF manipulators; and more recently,~\citet{kong2020fuzzy} integrate fuzzy-logic approximation with barrier Lyapunov functions to enforce semiglobal finite-time stability of $n-$link robots subject to actuator saturation and time-varying workspace constraints.
%\footnote{See~\cite{galicki2015finite, zhao2010robust}, and~\cite{kong2020fuzzy} for examples of this approach.}.
Consequently, there has been significant interest in verifying finite-time stability and establishing bounds on the settling-time function $T(x)$.

Methods for establishing finite-time stability are typically based on the use of Lyapunov functions and the comparison principle -- an approach established in the work of~\cite{bhat2000finite}. %\footnote{\cite{bhat2000finite,bhat1998continuous, haimo1986finite, venkataraman1990terminal, moulay2006finite}}.
%Methods for establishing finite-time stability typically proceed by constructing a Lyapunov function whose derivative along the system trajectories satisfies a negative-power inequality and then applying the comparison principle. In particular, ~\citep{bhat2000finite,bhat1998continuous, moulay2006finite} suppose the existence of a continuous, positive definite $V$ such that $\dot{V} \leq -\mu V^\gamma, \quad \mu >0, \gamma \in (0,1),$ which guarantees that the origin is finite-time stable.
Because finite-time stability of an equilibrium necessarily implies a non-Lipschitz vector field, however, the Lyapunov framework for finite-time stability is more nuanced than in the case of asymptotic stability notions and behaviour of a Lyapunov function near the equilibrium results in multiple Lyapunov characterizations (e.g. Lyapunov conditions in~\citet{bhat2000finite} differ significantly from those in~\cite{roxin1966finite}).

Moreover, although the Lyapunov framework for finite-time stability is relatively well-established, and has been used in an ad hoc manner to study consensus~\citep{wang2010finite}, barrier certificates~\citep{srinivasan2018control}, and twisting  algorithms~\citep{mendoza2017idea} et c., this Lyapunov framework has not resulted in model-agnostic polynomial-time algorithms for testing finite-time stability in the same way that Sum of Squares (SoS) is used to test rational stability of polynomial vector fields. The goal, then, is to establish a method for using SoS to analyze finite-time stable systems without restriction on the structure of the Lyapunov function.

%Moreover, although the Lyapunov framework for finite-time stability is relatively well-established, and has been used to study consensus~\citep{wang2010finite}, barrier certificates~\citep{srinivasan2018control}, and twisting  algorithms~\citep{mendoza2017idea} et c., and nonsmooth gradient flows in network consensus~\citep{cortes2006finite} -- it has not yet yielded polynomial-time algorithms for testing finite-time stability in the same way that Sum of Squares (SoS) is used to test rational stability of polynomial vector fields.

%Although Polya's theorem-based algorithms for finite-time stability analysis~\citep{sanchez2018sos, sanchez2014constructive} offer an SoS alternative, they impose a homogenizing structure on the Lyapunov function.

%In contrast, for finite-time stable homogeneous systems,a SoS method for constructing Lyapunov functions -- leveraging homogeneity properties through the application of P\'olya's theorem -- has been establsihed in \citep{efimov2018homogeneous,sanchez2014constructive,sanchez2019design}. The goal of this paper, then, is to establish such an algorithm to systematically certify finite-time stability.

For systems governed by polynomial vector fields, SoS programming, in conjunction with converse Lyapunov theorems, has proven to be an effective tool for stability certification \citep{topcu2008local, jarvis2005control}. Unfortunately, however, SoS methods are typically limited to analysis of systems with polynomial or rational vector fields and finite-time stable systems are typically characterized by non-polynomial dynamics. Attempts to extend SoS to non-polynomial vector fields include methods such as recasting, wherein the solutions of a non-polynomial system are embedded in a larger state-space with polynomial dynamics~\citep{savageau1987recasting} -- to which SoS methods can be applied~\citep{papachristodoulou2005analysis}. However, such methods are conservative in that instability of the recast system does not imply instability of the original system. More significantly, of course, polynomial vector fields cannot be finite-time stable and hence such methods cannot be applied to the construction of finite-time stability tests using SoS.

Recently, in~\citep{efimov2018homogeneous} and~\citep{sanchez2019design}, some progress has been made in this direction through the use of alternatives to SoS such as Polya's condition for positivity of homogeneous forms. In these works, either homogeneity is imposed on the vector field or the structure of the Lyapunov function is chosen to be homogenizing -- implying either a restriction of the class of systems or conservatism in the resulting stability conditions. By contrast, the methods presented here do not impose restrictions on the vector field or Lyapunov function, but instead use a state transformation to render the resulting stability conditions amenable to computation using SoS. This approach also provides accurate bounds on the settling time function and region of attraction on which that settling time function is valid.

%on SoS alternatives such as Polya's theorem exist for finite-time stable systems\footnote{See~\citep{efimov2018homogeneous} and \citep{sanchez2019design}}, these approaches impose a certain homogenizing structure on the Lyapunov function.

Specifically, we consider the case of a vector field defined by non-polynomial terms (e.g. fractional powers and signum functions). Then, instead of homogenizing or recasting the dynamics using a polynomial vector field, we take an alternative approach wherein we first pose the problem of finding a Lyapunov function which verifies finite-time stability -- Sec.~\ref{sec:problem}. This formulation includes inequalities with non-polynomial terms. However, as demonstrated in Sec.~\ref{sec:results}, a change of variables in the Lyapunov inequalities allows the non-polynomial inequality constraints to be represented by polynomial inequalities and enforced using SoS. A solution to these polynomial inequalities is then used to obtain a non-polynomial Lyapunov function which verifies the original stability conditions aned bounds the settling time-- Sec.~\ref{sec:finite_SOS_test}. This approach allows us to avoid the conservatism imposed by homogenization or recasting of the non-polynomial dynamics and results in accurate tests for finite-time stability and accurate bounds on the associated settling-time function, $T(x)$. In Sec.~\ref{sec:numerical} we verify the results using several test cases and use numerical simulation to evaluate the resulting bounds on settling time.
%%%%%%%%%%%%%%%%%%%%%%%%%%%%%%%%%%%%%%%%%%%%%%%%%%%%%%%%%%%%%%%%%%%%%%%%%%%%%%%%%%%%%%%%%%%%%%%%%%%%%%%%%%%%%%%%%%%%%%%%%%%%%%%%%%%%%%%%%%%%%%%%%%%%%%%%%%%%%%%%%%%%%%%%%%%%%%%%%%%%%%%%%%%%%%%%%%%%%%%%%%%%%%%%%%%%%%%%%%%%%%%%%
\section{Notation}\label{sec: notation}\vspace{-2mm}
$\mathbb{R}^n$, $\mathbb{R}_+^n$, and $\mathbb{N}^n$ denote the space of $n$-dimensional vectors of real, positive real, and natural numbers, respectively. $\mathbb{R}[x]$ as the set of real-valued polynomials in variables $x$.
A neighborhood, $\mcl N$ of a point, $x$ is any set which contains an open set $U$ such that $x \in U \subseteq \mcl N$.
In Thm.~\ref{higherdimvarsub}, we introduce elementwise vector operations for multiplication, $\sign$ and exponent. Specifically, for $x \in \mathbb{R}^n$ and $q \in \mathbb{N}^n$, we define the vector $y=x^q$ as $y_i=x_i^{q_i}$. We also define vector-valued element-wise $\sign$ and absolute value functions where if $x\in \R^n$ and $y=\sign(x)$, then $y_i=sign(x_i)$ and if $y=|x|$ then $y_i=|x_i|$. Finally, we define the element-wise (Hadamard) product of vectors so that if $y=z \cdot x$, then $y_i=z_ix_i$.
%%%%%%%%%%%%%%%%%%%%%%%%%%%%%%%%%%%%%%%%%%%%%%%%%%%%%%%%%%%%%%%%%%%%%%%%%%%%%%%%%%%%%%%%%%%%%%%%%%%%%%%%%%%%%%%%%%%%%%%%%%%%%%%%%%%%%%%%%%%%%%%%%%%%%%%%%%%%%%%%%%%%%%%%%%%%%%%%%%%%%%%%%%%%%%%%%%%%%%%%%%%%%%%%%%%%%%%%
\section{Sum of Squares Decomposition}\label{sec:SOS}\vspace{-2mm}
In this section, we provide a brief overview of SoS polynomials and how  semidefinite programming can be utilized to verify the existence of SoS decomposition. Specifically, polynomial $p(x)\ge 0$ for all $x \in\R^n$ if $p$ is Sum-of-Squares.

\begin{defn}
	A polynomial $ p(x) \in \mathbb{R}[x] $ is called a \emph{sum of squares (SoS)}, denoted $ p \in \Sigma_s $, if there exist polynomials $h_i(x) \in \mathbb{R}[x] $  such that $	p(x) = \sum_{i} h_i(x)^2$.
\end{defn} 
The existence of a SoS representation of a polynomial, $p$, can be tested using semidefinite programming~\citep{parrilo2000structured}. Specifically, $p(x)$ of degree $2d$ is SoS if and only if there exists a positive semidefinite matrix $Q$ such that $p=Z_d(x)^T Q Z_d(x)$ where $Z_d(x)$ is the vector of monomials in $x$ of degree $d$ or less. In  the case where the polynomial inequality is required to hold on a semi-algebraic set (i.e. $p(x)\ge 0$ for all $x\in S$), we may use Positivstellensatz results. Specifically, a set $S $ is called \emph{semi-algebraic} if it can be represented using polynomial equality and inequality constraints as
	\[
	S = \left\{ x :
	\begin{aligned}
		& g_i(x) \ge 0, \quad \forall i = 1, \ldots, k, \quad k \in \mathbb{N}, \\
		& h_j(x) = 0, \quad \forall j= 1, \ldots, m, \quad m \in \mathbb{N}
	\end{aligned}
	\right\}.
	\]

\begin{lem}
	Given $ S = \{x : g_i(x) \geq 0, h_j(x) = 0 \} $, suppose that there exist $s_i \in \Sigma_s $ and $ t_j \in \mathbb{R}[x] $ such that
	\[
	p(x) = s_0(x) + \sum^k_{i=1} s_i(x) g_i(x) + \sum^m_{j=1} t_j(x) h_j(x),
	\]
	where $ k, m \in \mathbb{N} $. Then $p(x)\ge 0$ for all $x \in S$.
\end{lem}\vspace{-2mm}
Necessity of this Positivstellensatz test in the case of strictly positive polynomials holds if $s_i,h_i$ satisify an additional precompactness condition~\citep[Lemma 4.1]{putinar1993positive}. All semi-algebraic sets used in this paper satisfy this condition. 
%%%%%%%%%%%%%%%%%%%%%%%%%%%%%%%%%%%%%%%%%%%%%%%%%%%%%%%%%%%%%%%%%%%%%%%%%%%%%%%%%%%%%%%%%%%%%%%%%%%%%%%%%%%%%%%%%%%%%%%%%%%%%%%%
\section{Problem Formulation}\label{sec:problem} \vspace{-2mm}

We consider nonlinear differential equations of the form
\begin{equation}\label{eq1}
	\frac{d}{d t}x(t) = f(x(t)), \quad x(0)=x_0,t \in [0,T), \quad T \in \mathbb{R}_+,
\end{equation}
where $x(t) \in \mathcal{D} \subseteq \mathbb{R}^n$ with $\mathcal{D}$ containing an open neighborhood of the origin. Furthermore, we suppose that $f: \mathcal{D} \mapsto \mathbb{R}^n$ is continuous  on $\mathcal{D}$ with $f(0) = 0$. However, we do not require $f$ to be locally Lipschitz at the origin. For this reason, we require existence and uniqueness of solutions except possibly at the origin, where such uniqueness is as defined in~\citep{bhat2000finite}.
In this case, we denote by $\phi_f :[0,T)  \times \mathcal{D} \backslash \{0\} \mapsto \mathcal{D} \backslash \{0\}$ the corresponding solution map for which
\[
\begin{aligned}
	&\frac{d}{dt} \phi_f (t,x) = f(\phi_f (t,x)), \quad &&\forall x \in \mathcal{D} \backslash \{0\} , t \in [0,T) , \\
	&\phi_f (0, x) = x,  \quad &&\forall x \in \mathcal{D} \backslash \{0\}.
\end{aligned}
\]
For a nonlinear system with the associated solution map, we may define a notion of finite-time stability using the framework proposed in~\citep{bhat2000finite}.
\begin{defn}\label{FTS}
	Given $f$ and $\mcl O\subset \mcl D$ containing an open neighborhood of the origin, we say the solution map $\phi_f :[0,T)  \times \mathcal{D}\backslash \{0\} \mapsto \mathcal{D} \backslash \{0\}$ is \emph{finite-time stable} on $\mcl O$ with settling-time function $T: \mathcal{O} \backslash \{0\} \rightarrow \mathbb{R}_+$ if it is stable in the sense of Lyapunov and
	for every $x \in \mathcal{O} \backslash \{0\}$, $\phi_f (t,x) \in \mathcal{O} \backslash \{0\}$ for all $t \in [0, T(x))$ and $\lim_{t \to T(x)} \phi_f(t,x) = 0.$
	If $\mathcal{O} = \mathcal{D} = \mathbb{R}^n$, we say that $\phi_f$ is \emph{globally finite-time stable}.
\end{defn}
While the above definition of finite-time stability allows for any valid settling-time function, to avoid ambiguity, henceforth we refer to \textbf{the} settling-time function as the function, $T$, which is defined as  $T(0)=0$ and
\(
T(x) = \inf \{ t \in \mathbb{R}_+ \mid \phi_f(t, x) = 0 \} .
\)

The goal, then, is to establish SoS-based conditions which prove finite-time stability and to obtain least upper bounds on the associated settling-time function. To achieve this goal, we rely on the Lyapunov characterization of finite-time stability and the associated settling-time function developed in~\cite[Theorem 4.2]{bhat2000finite}:

\begin{thm}\label{bglem}
	Suppose there exists a continous function $V: \mathcal{D} \mapsto \mathbb{R}$ such that the following conditions hold:
	(i) \, $V$ is positive definite.
		(ii) \, There exist real numbers $\mu>0$ and $\gamma \in (0,1)$ and an open neighborhood $\Omega \subseteq \mathcal{D}$ of the origin such that
		\begin{equation}
			\dot{V}(x) + \mu V(x)^\gamma \leq 0, \quad x \in \Omega \backslash \{0\}.
		\end{equation}
	Then the solution map  $\phi_f:[0,T)  \times \mathcal{D} \backslash \{0\} \rightarrow \mathcal{D} \backslash \{0\}$ of \eqref{eq1} is finite-time stable. Moreover, if $\mathcal{O}$ is contained in a sublevel set of $V$ contained in $\Omega$ and $T$ is the settling-time function, then
	\begin{equation}\label{lem12}
		T(x) \leq \frac{1}{\mu(1-\gamma)} V(x)^{1-\gamma}, \quad x \in \mathcal{O},
	\end{equation}
	and $T$ is continuous on $\mathcal{O}$. If in addition $\mathcal{D}=\mathbb{R}^n, V$ is proper, and $\dot{V}$ takes negative values on $\mathbb{R}^n \backslash \{0\},$ then the origin is a globally finite-time stable equilibrium of \eqref{eq1}.
\end{thm}
While Thm.~\ref{bglem} provides sufficient conditions for finite-time stability in terms of a continuous Lyapunov function, these conditions are expressed using a fractional exponent, $\gamma \in (0,1)$. %Furthermore, finite-time stable systems are defined by non-polynomial vector fields -- often containing fractional exponents. 
As such, typical SoS and polynomial programming methods cannot be applied. In the following section, we show how to reformulate the conditions of Thm.~\ref{bglem} using polynomials and polynomial inequalities.

%%%%%%%%%%%%%%%%%%%%%%%%%%%%%%%%%%%%%%%%%%%%%%%%%%%%%%%%%%%%%%%%%%%%%%%%%%%%%%%%%%%%%%%%%%%%%%%%%%%%%%%%%%%%%%%%%%%%%%%%%%%%%%%%
\section{A Polynomial Reformulation of the Finite-Time Stability Problem}\label{sec:results} \vspace{-2mm}
The Lyapunov conditions in Thm.~\ref{bglem} contain fractional exponents (e.g. $\dot V \le -V^{1/2}$) and finite-time stable systems typically include vector fields with fractional exponents (e.g. $\dot x=-x^{1/5}$). Such fractional terms prevent a straightforward application of polynomial optimization and SoS. In this section, we show how the conditions of Thm.~\ref{bglem} can be enforced without the use of fractional exponents.

To keep the exposition clear, we first state and prove the result for the scalar case, where $f : \R \rightarrow \R$.  Building on the scalar case, we provide a similar result for the multivariate case in Thm.~\ref{higherdimvarsub}.
%The proof for multivariate vector fields in Thm.~\ref{higherdimvarsub} is similar, but requires additional terms and notation.

\begin{thm}\label{varsub}
Let $\mathcal{D} \subseteq \mathbb{R}$ and $f:\mathcal{D} \mapsto \mathbb{R}$ be continuous. Suppose for $\Omega \subset \R$, containing an open neighborhood of the origin, there exist a continuously differentiable $V:\Omega \rightarrow \mathbb{R}$, positive constants $\mu, k$,  and integers $p,q, r \in \N$ such that $V(0)=0, V(x)>0$ for any $x \in \Omega \backslash \{0\}$,  $V(x)\le k |x|^r $ for $x \in \Omega$ ,  and
\begin{equation}\label{conv1}
		\nabla_x V(x) f(\sign(x)|x|^q)  \leq -\mu  |x|^p, \quad \forall x \in \Omega .
\end{equation}
%%%%%%%%%%%%%%%%%%%%%%%%%%%%%%%%%
Let $\tilde{\Omega } := \{z \in \mathcal{D}: \sign(z)|z|^{1/q} \in \Omega \}$, $\gamma =\frac{1+p-q}{r} \in \mathbb{R}_+$ and $\tilde \mu = \frac{\mu}{q k^\gamma}$.
	Then there exists  a continuous $\tilde V:\tilde{\Omega} \rightarrow \mathbb{R}$, continuously differentiable on $\tilde{\Omega} \backslash \{0\},$ such that $\tilde V(0)=0, \tilde V(z)>0$ for  any $z \in \tilde \Omega \backslash \{0\}$, and
	\begin{equation}\label{conv2}
		\nabla_z \tilde{V}(z)  f(z) \le - \tilde \mu \tilde{V}(z)^\gamma , \quad \forall z \in \tilde{\Omega} \backslash \{0\}.
	\end{equation}
\end{thm}
\vspace{-2mm}
\begin{pf}
	For $z \in \tilde{\Omega}$, define
	\(
	\tilde{V}(z) = V(\sign(z)|z|^{1/q}).
	\)
	First, $\tilde{V}(0) = V(0) = 0$. Next, for any $z\in \tilde \Omega,$ $\sign(z)|z|^{1/q} \in \Omega$ and hence $\tilde V $ is well-defined on $\tilde \Omega$. Furthermore,  $\tilde{V}(z) = V(\sign(z)|z|^{1/q}) > 0$ for all $z \in \tilde{\Omega}.$
	Next, since  $V$ is continuously differentiable on $\Omega$ and
\[
	\frac{\partial}{\partial z}\sign(z)|z|^{1/q} = \frac{\partial}{\partial z} z|z|^{1/q-1} =\frac{1}{q}|z|^{1/q-1}
\]
is continuous on $\tilde \Omega \backslash \{0\}$, we have $\tilde{V}(z) = V(\sign(z)|z|^{1/q})$ is continuous on $\Omega$ and continuously differentiable on $\tilde \Omega \backslash \{0\}$.

	Finally, using the chain rule and the expression for 
	\(\frac{\partial}{\partial z}\sign(z)|z|^{1/q} \), we have:
	\[
	\nabla_z \tilde{V}(z) = \nabla_x V(\sign(z)|z|^{1/q}) \frac{1}{q}|z|^{1/q-1} .
	\]
	We now observe that
	\(
	\sign(\sign(z)|z|^{1/q})|\sign(z)|z|^{1/q}|^q=\sign(z)|z|=z
	\)
	and hence for $z \in \Omega$,~\eqref{conv1} implies
	\[
	\nabla_x V(x)\vert_{x=\sign(z)|z|^{1/q}}  f(z)  \leq -\mu  |\sign(z)|z|^{1/q}|^p= -\mu  |z|^{p/q}.
	\]
	Thus
	\begin{align*}
		\nabla_z \tilde{V}(z) f(z)
		&= \nabla_x V(x)\vert_{x=\sign(z)|z|^{1/q}} \frac{1}{q}|z|^{1/q-1} f(z) \\
		&\le -\mu \frac{1}{q}|z|^{1/q-1} |z|^{p/q} = -\frac{\mu}{q}|z|^{\frac{1+p-q}{q}}.
	\end{align*}
	
	However, we know that $V(x) \leq  k |x|^r$ with $r = \frac{1+p-q}{\gamma}$ which implies
	\begin{align*}
	\tilde V(z) &= V(\sign(z)|z|^{1/q}) \\ &\leq  k  |\sign(z)|z|^{1/q}|^{\frac{1+p-q}{\gamma}} =k |z|^{\frac{1+p-q}{q\gamma}}.
	\end{align*}
	Rasing both sides to the power of  $\gamma$, we have:
	 \(
	 \tilde V(z)^\gamma \leq k^\gamma |z|^{\frac{1+p-q}{q}}.
	 \)
	Hence, for any $z \in \tilde \Omega \backslash \{0\}$ we obtain:
	\[
	\nabla_z \tilde{V}(z) f(z) \le -\frac{\mu}{q}|z|^{\frac{1+p-q}{q}}\le -\frac{\mu}{qk^\gamma}  \tilde V(z)^\gamma
	\]
	as desired. \hfill $\square$ 
\end{pf} 
\begin{note}
In the formulation of Thm.~\ref{varsub}, $x\in \Omega$ is the varaiable used in the polynomial inequalities which will be enforces using SoS while variable $z\in \tilde{\Omega}$ is the state of the original system and $\tilde V(z)$ is the resulting Lyapunov function for that system.
\end{note}

Theorem~\ref{varsub} provides alternative conditions under which Theorem~\ref{bglem} can be used to prove finite-time stability and bound the settling-time function. If the vector field, $f$ contains fractional exponents, e.g. $f(x)=-\sign(x)|x|^{1/2}-x^{1/3}$ then by choosing $q$ to be the least common denominator of these fractional terms, $f(\sign(z)|z|^q)$ the fractional terms are eliminated -- e.g. $f(\sign(z)|z|^6)=-x^{3}-x|x|$. While these conditions are not entirely polynomial, they can then be tested using polynomial optimization as described in Sec.~\ref{sec:finite_SOS_test} and illustrated in Sec.~\ref{sec:numerical} (Example~\ref{ex:1}).% For the multivariate case, however, the conditions are slightly more involved, as seen in the following subsection.

%%%%%%%%%%%%%%%%%%%%%%%%%%%%%%%%%%%%%%%%%%%%%%%%%%%%%%%%%%%%%%%%%%%%%%%%%%%%%%%%%%%%%%%%%%%%%%%%%%%%%%%%%%%%%%%%%%%%%%%%%%%%%%%%%%%%%%%%%%%%%%%%%%%%%%%%%%%%%%%%%%%%%%%%%%%%%%%%%%%%%%%%%%%%%%%%%%%%%%%%%%%%%%%%%%%%%%%%%%%%%%%%%%%%%%

\subsection{Finite-Time Stability for Multivariate Systems }\label{sec:results2}\vspace{-2mm}
In Theorem~\ref{varsub} we have shown, in the scalar case, how fractional terms in the Lyapunov test for finite-time stability may be eliminated by a substitution $x \mapsto \sign(z)|z|^q$.
In the multivariate case, we now present similar conditions, although in this case, we allow for multiple variable mappings -- e.g. $x_1 \mapsto \sign(z_1)|z_1|^2$ and $x_2 \mapsto \sign(z_2)|z_2|^6$. As a result, however, there is an additional step in modifying the vector field -- i.e. the conditions are expressed in terms of a modified vector field $\tilde f_i(z):=\frac{1}{q_i} f_i(z)|z_i|^{1/q_i-1}$.
%Having stated the main result in the scalar case, we now provide a multivariate version, at the cost of some additional notation (See Sec.~\ref{sec: notation}). %The corollary presented herein constitutes a generalization of Theorem~\ref{varsub} from the scalar setting to the multivariate case, extending the methodology to systems defined over $\mathbb{R}^n.$ The key objective remains the reformulation of fractional Lyapunov inequalities in Theorem~\ref{bglem}  as a polynomial optimization program.

\begin{thm}\label{higherdimvarsub}
	Let $\mathcal{D} \subseteq \mathbb{R}^n$ and $f: \mathcal{D} \mapsto \mathbb{R}^n$ be continuous.
	Suppose $\Omega \subset \mathbb{R}^n$ contains an open neighborhood of the origin and there exist a continuously differentiable
	$V: \Omega  \rightarrow \mathbb{R}$, a vector $q \in \mathbb{N}^n$, and scalars $p, \mu, k,r \in \mathbb{R}_+$
	such that $V(0) = 0$, $V(x) > 0$ for any $x \in \Omega \backslash \{0\}$,
	and $V(x) \leq  k\norm{x}^{r}$ for $x \in \Omega$.
	Moreover, suppose
	\begin{equation}\label{conv1_higher}
		\nabla_x V(x)^T \tilde f(\sign(x) \cdot |x|^q) \leq  - \mu \norm{x}^{p}, \quad \forall x \in \Omega,
	\end{equation}
	where $\tilde f_i(z):=\frac{1}{q_i} f_i(z)|z_i|^{1/q_i-1}$.
	Let $\tilde{\Omega} = \{z \in \mathcal{D}: \sign(z) \cdot |z|^{1/q} \in \Omega \}$, $\gamma = \frac{p}{r}$, and $\tilde{\mu} = \frac{\mu}{k}$.
	Then there exists a continuous $\tilde{V}:\tilde{\Omega} \rightarrow \mathbb{R}$, continuously differentiable on $\tilde{\Omega} \backslash \{0\},$ such that
	$\tilde{V}(0)=0$,
	$\tilde{V}(z) > 0$ for any $z \in \tilde \Omega \backslash \{ 0\}$ and
	\begin{equation}\label{conv2_higher}
		\nabla_z \tilde{V}(z)^T f(z) \le  - \tilde{\mu} \tilde{V}(z)^{\gamma} \quad \forall z \in \tilde{\Omega} \backslash \{0\}.
	\end{equation}
\end{thm}
\begin{pf}
	%Let $V:\Omega \rightarrow \mathbb{R}$ be continuously differentiable and suppose that there exist a vector $q \in \mathbb{N}^n$, and scalars $p, \mu, k,r \in \mathbb{R}_+$ such that the conditions of the theorem statement hold.
	For $z \in \tilde{\Omega}$, define
	\(
	\tilde{V}(z) = V(\sign(z) \cdot |z|^{1/q}).
	\)
	%where $\sign(z) \cdot |z|^{1/q} = [\sign(z_1) |z_1|^{1/q_1}; \cdots; \sign(z_n) |z_n|^{1/q_n} ]$.
	First, $\tilde{V}(0) = V(0) = 0$. Next, for any $z\in \tilde \Omega$, $\sign(z)\cdot |z|^{1/q} \in \Omega$ and hence $\tilde V $ is well-defined on $\tilde \Omega$. Furthermore, $V(x)>0$ for all $x \in {\Omega} \backslash \{0\}$ implies $\tilde{V}(z) = V(\sign(z) \cdot |z|^{1/q}) > 0$ for all $z \in \tilde{\Omega}\backslash \{0\}$.
	Next, since $\tilde{V}(z) = V(\sign(z)\cdot |z|^{1/q})$, $V$ is continuously differentiable on $\Omega$ and, as in the scalar case,
\(	
\frac{\partial}{\partial z_i}\sign(z_i)|z_i|^{1/q_i} =\frac{1}{q_i}|z_i|^{1/q_i-1}
\)
is continuous on $\tilde \Omega \backslash \{0\}$, we have that $\tilde{V}(z) = V(\sign(z) \cdot |z|^{1/q})$ is continuous on $\tilde{\Omega}$ and continuously differentiable on $\tilde \Omega \backslash \{0\}$.
Finally, using the chain rule and the expression for \(\frac{\partial}{\partial z}\sign(z)|z|^{1/q} \), we have:
	\[
	\frac{\partial \tilde{V}(z)}{\partial z_i} = \frac{\partial V(x)}{\partial x_i}\vert_{x_i=\sign(z_i) |z_i|^{1/q_i}}  \frac{1}{q_i} |z_i|^{1/q_i-1}.
	\]
	Hence $\nabla_z \tilde{V}(z)^T f(z)  = \sum_{i} \frac{\partial \tilde{V}(z)}{\partial z_i} f_i(z)$ implies
\[		\nabla_z \tilde{V}(z)^T f(z) % & = \sum_{i=1}^{n} \frac{\partial \tilde{V}(z)}{\partial z_i} f_i(z) \\
		=\sum_{i=1}^{n} \partial_i V(\sign(z) \cdot |z|^{1/q})f_i(z)\frac{1}{q_i} |z_i|^{1/q_i-1}
		%&=\sum_{i=1}^{n} \partial_i V(\sign(z) \cdot |z|^{1/q})\tilde f_i(z).
\]  %For each $i$ define  $\tilde f_i(z):=\frac{1}{q_i} f_i(z)|z_i|^{1/q_i-1}.$
where $\partial_i V$ denotes derivative with respect to the $i^{th}$ argument. Now by assumption we have for $x \in \Omega$:
	\[
	\sum_{i=1}^{n} \partial_i V(x)\tilde f_i(\sign(x) \cdot |x|^{q}) \leq  -\mu \norm{x}^{p}
	\]
and since $	\sign(\sign(z_j)|z_j|^{1/q_j})|\sign(z_j)  |z_j|^{1/q_j}|^{q_j}=z_j$ 
	 for $z \in \tilde \Omega$, we have $\norm{\sign(z) \cdot |z|^{1/q}}=\norm{|z|^{1/q}}$ implies
%	\begin{align*}
%		\sum_{i=1}^{n} \partial_i V(\sign(z) \cdot |z|^{1/q})\tilde f_i(z) &\leq  -\mu  \norm{\sign(z_i) |z_i|^{1/q_i}}^{p}, \\
%		&=-\mu \norm{|z_i|^{1/q_i}}^p.
%	\end{align*}
\[
		\sum_{i=1}^{n} \frac{ \partial V}{\partial_{z_i}}(\sign(z) \cdot |z|^{1/q})\tilde f_i(z) \leq -\mu \norm{|z|^{1/q}}^p
\]	where recall from notation that vector $|z|^{1/q}$ has elements $|z_i|^{1/q_i}$. Hence we conclude that
\[
		\nabla_z \tilde{V}(z)^T f(z)=\sum_{i=1}^{n} \partial_i V(\sign(z) \cdot |z|^{1/q}) \tilde f_i(z)\le-\mu\norm{|z|^{1/q}}^p.
	\]
    However, we know that for $x \in \Omega$:
	\(
	V(x) \leq k\norm{x}^{r}.
	\)
	Now for any $z\in \tilde \Omega$, there exists $x \in \Omega$ such that $x_i = \sign(z_i)|z_i|^{1/q_i}$, hence
\[
		\tilde{V}(z)= V(\sign(z) \cdot |z|^{1/q}) \le k \norm{|z|^{1/q}}^r.
\]
	Raising both sides to the power $\gamma=p/r$ and dividing by $k$, we obtain
	\(
	\frac{1}{k}\tilde{V}(z)^{\gamma}\le \norm{|z|^{1/q}}^p
	\)
	which implies
	\(
	-\norm{|z|^{1/q}}^p\le - \frac{1}{k}\tilde{V}(z)^{\gamma}.
	\)
	Therefore,  for $ z \in \tilde{\Omega} \backslash \{0\}$, we conclude:
	\[
	\sum_{i=1}^{n} \frac{\partial \tilde{V}(z)}{\partial z_i} f_i(z)\le-\mu\norm{|z|^{1/q}}^p\le - \frac{\mu}{k}\tilde{V}(z)^{\gamma}
	\]
	which implies
	\(
	\nabla_z \tilde{V}(z)^T f(z) \le  - \tilde{\mu} \tilde{V}(z)^{\gamma}
	\)
	as desired.  \hfill $\square$ 
\end{pf} 
%%%%%%%%%%%%%%%%%%%%%%%%%%%%%%%%%%%%%%%%%%%%%%%%%%%%%%%%%%%%%%%%%%%%%%%%%%%%%%%%%%%%%%%%%%%%%%%%%%%%%%%%%%%%%%%%%
Similar to Thm.~\ref{varsub}, Thm.~\ref{higherdimvarsub} provides alternative conditions under which the stability conditions of Thm.~\ref{bglem} are satisfied. Like in the scalar case, through judicious choice of $q_i$, the conditions of Thm.~\ref{higherdimvarsub} eliminate fractional terms from the Lyapunov conditions. Unlike Thm.~\ref{varsub}, however, the $\tilde f(\sign(x)\cdot |x|^q)$ function in the conditions of Thm.~\ref{higherdimvarsub} may contain rational terms. For example, if $f_1(x,y)=-\sign(x)\sqrt{|x|}-\sqrt[3]{y}$ and we choose $q_1=2$ and $q_2=3$, then
$\tilde{f}_1(\sign(z)\cdot |z|^q)=-\frac{1}{2}\frac{z_1+z_2}{|z_1|}$. Fortunately, these rational conditions can also be tested using polynomial optimization as described in Section~\ref{sec:finite_SOS_test} and illustrated in Section~\ref{sec:numerical} (Example~\ref{ex:2}).

Having formulated alternative Lyapunov stability conditions, we now combine Thm.~\ref{higherdimvarsub} with Thm.~\ref{bglem} to formally show that these imply that the solution map, $\phi_f$, associated with the vector field, $f$, is finite-time stable and provides a bound on the associated settling-time function.

% We now that reformulate the Lyapunov theorem for finite time stability, Theorem~\ref{bglem} as polynomial Lyapunov inequalities, we can conclude that origin is a finite time stable equilbirum of \eqref{eq1} with a polynomial settling time bound.
\begin{cor}\label{Cor01}
	Let $\mathcal{D} \subseteq \mathbb{R}^n$ and $f: \mathcal{D} \mapsto \mathbb{R}^n$ be continuous. Suppose for $\Omega \subset \R^n$,  containing an open neighborhood of the origin, there exist a continuously differentiable $V:\Omega \rightarrow \mathbb{R}$, a vector $q \in \mathbb{N}^n$, and scalars  $k, p,r, \mu \in \mathbb{R}_+$ such that $V(0)=0, V(x)>0$ for any $x \in \Omega \backslash \{0\}$,  $V(x)\le k \|x\|^{r}$ for $x \in \Omega$, and
	\begin{equation}
		\nabla_x V(x)^T \tilde f(\sign(x) \cdot |x|^q)  \leq -\mu  \|x\|^p, \quad \forall x \in \Omega,
	\end{equation}
	where $\tilde f_i(z):=\frac{1}{q_i} f_i(z)|z_i|^{1/q_i-1}$.
	Let $\gamma = \frac{p}{r}$, $\tilde{\mu} = \frac{\mu}{k}$, $\tilde{\Omega} = \{z \in \mathcal{D}: \sign(z)\cdot |z|^{1/q} \in \Omega \} $, and $\tilde{V}(z) := V(\sign(z)\cdot |z|^{1/q})$.
Then the solution map $\phi_f:[0,T)  \times \mathcal{D} \backslash \{0\} \rightarrow \mathcal{D} \backslash \{0\}$ of \eqref{eq1} is finite-time stable. Additionally, if $\mathcal{O}$ is contained in a sublevel set of $\tilde V$ that is itself contained in $ \tilde \Omega$, and $T$ is the settling-time function, then for $z \in \mathcal{O}$:
\begin{equation}\label{Cor}
	T(z) \leq \frac{1}{\tilde  \mu(1-\gamma)} \tilde V(z)^{1-\gamma}.
	\end{equation}
\end{cor}
\begin{pf}
Suppose that the conditions of the corollary statement are satisfied. Then by Thm.~\ref{higherdimvarsub}, if $\tilde{\Omega } = \{z \in \mcl D: \sign(z)\cdot |z|^{1/q} \in \Omega \}$, $\gamma =\frac{p}{r}$ and $\tilde \mu =\frac{\mu}{k}$, there exists a continuous $\tilde{V}:\tilde{\Omega} \rightarrow \mathbb{R}$, continuously differentiable on $\tilde{\Omega} \backslash \{0\},$ such that $\tilde{V}(0)=0$, $\tilde{V}(z) > 0$ for any $z\in \Omega \backslash \{0\}$ and
$\nabla_z \tilde{V}(z)^T  f(z) \le  - \tilde \mu \tilde{V}(z)^{\gamma}$ for all $z \in \tilde{\Omega} \backslash \{0\}$. Since $\gamma \in (0,1)$ the conditions of Thm.~\ref{bglem} are satisfied and therefore the solution map $\phi_f$ of \eqref{eq1} is finite-time stable , with the settling-time function $T$ that satisfies~\eqref{Cor}.  \hfill $\square$ 
\end{pf} 
\section{SoS Conditions for Finite-Time Stability}\label{sec:finite_SOS_test}\vspace{-2mm}
%%%%%%%%%%%%%%%%%%%%%%%%%%%%%%%%%%%%%%%%%%%%%%%%%%%%%%%%%%%%%%%%%%%%%%%%%%%%%%%%%%%%%%%%%%%%%%%%%%%%%%%%%%%%%%%%%%%%%%%%%%%%%%%%%%%%%%%%%%%%%%%%%%%%%%%%%%%%%%%%%%%%%%%%%%%%%%%%%%%%%%%%%%%%%%%%%%%%%%%%%%%%%%%%%%%%%%%%%%%
Thm.~\ref{varsub} and Thm.~\ref{higherdimvarsub} were motivated by a desire to use polynomial programming to prove finite-time stability and bound the settling-time function. In the following theorem, we propose such a polynomial programming problem, based on the use of SoS and Positivstellensatz results to enforce the conditions of Thm.~\ref{higherdimvarsub} on some compact semi-algebraic set containing the origin.

%We now reformulate the polynomial Lyapunov inequalities in Corollary~\ref{Cor01} as SoS, thereby enabling the use of semidefinite programming to  compute a Lyapunov function $V(x)$ that satisfies the proposed sufficient conditions  for verifying finite-time stability and providing precise bounds on the settling time.  The following theorem formalizes this idea, showing that if certain conditions are met, including the existence of a valid inequality on a semi-algebraic set, then the system exhibits finite-time stability. Moreover, it provides an explicit way to compute bounds on the settling time using SoS optimization.

\begin{prop}\label{SoSApp}
	Let $\mathcal{D} \subseteq \mathbb{R}^n$ and $f:\mathcal{D} \mapsto \mathbb{R}^n$ be continuous. Suppose there exist a polynomial $V: \Omega \mapsto \mathbb{R}$, vectors $q, \lambda \in \mathbb{N}^n$, integers $\tau, d,p \in \mathbb{N}$, scalars $\mu, \epsilon, k \in \mathbb{R}_+$, SoS $s_{i}, t_{i} $, and some semi-algebraic set $\Omega \subset \mathbb{R}^n$ containing the origin with valid inequalities $g_{i}$ such that
$ V(0) = 0$ , $V(x) - \epsilon \|x \|^{2 \tau} \in \Sigma_s$, and 
\begin{align}
&		k \| x \|^{2d} - V(x) - \sum\nolimits_{i=1}^n s_{i}(x) g_{i}(x) \in \Sigma_s,  \label{eq:cond3}\\
		&-\nabla_x V(x)^T \tilde f (\sign(x)\cdot |x|^{q}) \prod_{i=1}^n |x_i|^{\lambda_i} - \mu \|x\|^{p} \prod_{i=1}^n |x_i|^{\lambda_i} \notag\\[-2mm]
			&\qquad \qquad\qquad \qquad\qquad\qquad- \sum_{i=1}^n t_{i}(x) g_{i}(x) \in \Sigma_s, \label{eq:cond4}
\end{align}
where $\tilde f_i(z):=\frac{1}{q_i} f_i(z)|z_i|^{1/q_i-1}$.
	Let  $\tilde{\mu} = \frac{\mu}{k}$, $\gamma=\frac{p}{2d}$, $\tilde{\Omega} = \{z \in \mathcal{D}: \sign(z)\cdot |z|^{1/q} \in \Omega \}$, and  $\tilde{V}(z) := V(\sign(z)\cdot |z|^{1/q})$.
	Then for $2d > p$, the solution map $\phi_f:[0,T)  \times \mathcal{D} \backslash \{0\} \rightarrow \mathcal{D} \backslash \{0\}$ of~\eqref{eq1} is finite-time stable. Furthermore, if $\mathcal{O}$  belongs to a  sublevel set of $ \tilde V$ that is contained in $\tilde \Omega$, and $T$ is the settling-time function, then
	\begin{equation}\label{cor01}
		T(z) \leq \frac{1}{\tilde \mu(1-\gamma)} \tilde V(z)^{1-\gamma} , \quad z  \in \mathcal{O}.
	\end{equation}
\end{prop}
%%%%%%%%%%%%%%%%%%%%%%%%%%%%%%%%%%%%%%%%%%%%%%%%%%%%%%%%%%%%%%%
\begin{pf}
	Suppose the conditions of the theorem are satisfied.
	Then for $x \in \Omega$,  the condition $V(x) -\epsilon \|x\|^{2 \tau} \in \Sigma_s$  ensures $V(x) > 0$ for $x \in \Omega \backslash \{0\}$.
	Moreover, for $x \in \Omega$, the condition $k\|x\|^{2d}-V(x) -\sum s_i(x)g_i(x) \in \Sigma_s$ implies that for $x \in \Omega$:
	\(
	V(x) \leq k\|x\|^{2d} -\sum s_i(x)g_i(x) \leq k \|x\|^{2d} .
	\)
	This also implies $V(0)=0$. Now $2d > p$ implies that $\gamma=\frac{p}{2d} \in (0,1)$ and for $x \in \Omega$,
	\(
	V(x)  \leq k\|x\|^{2d} :=k\|x\|^{\frac{p}{\gamma}}.
	\)
	
	Since $\prod_{i=1}^n |x_i|^{ \lambda_i} \geq 0$, $s_i, t_i\in \Sigma_s$ and $g_i(x)\ge 0$ for $x \in \Omega$, Condition~\eqref{eq:cond4} implies that for $x \in \Omega,$
\[
	\nabla_x V(x)^T \tilde f(\sign(x) \cdot |x|^q) \prod_{i=1}^n |x_i|^{ \lambda_i} \leq - \mu \|x\|^p \prod_{i=1}^n |x_i|^{ \lambda_i},
\]
which implies
	\(
	\nabla_x V(x)^T \tilde f(\sign(x) \cdot |x|^q)  \leq -\mu  \|x\|^p
	\)
	for all $x \in \Omega$. Hence, the conditions of Cor.~\ref{Cor01} are satisfied for the given $k,\mu,q,p,\Omega$ and for $r=2d$. We conclude that if $\gamma = \frac{p}{r}$, $\tilde{\mu} = \frac{\mu}{k}$, $\tilde{V}(z) := V(\sign(z)\cdot |z|^{1/q})$, $\tilde{\Omega} = \{z \in \mathcal{D}: \sign(z) \cdot |z|^{1/q} \in \Omega \}$, and $\mathcal{O}$ belongs to a  sublevel set of $ \tilde V$ that is contained in $\tilde \Omega$ then the solution map $\phi_f$ of \eqref{eq1} is finite-time stable with the settling-time function as in the statement of the proposition. \hfill $\square$ 
\end{pf} 
\begin{note}
Condition~\eqref{eq:cond3} of Prop.~\ref{SoSApp} implies that the upper bound on the function $V$ only holds locally. This allows for $d < \tau$ since in such a case it is not possible for the upper bound to hold globally.
\end{note}

\begin{note}The conditions of Prop.~\ref{SoSApp} are formulated ostensibly as an SoS programming problem with polynomial variables $V,s_i,t_i$. Furthermore, If the $q_i$ are chosen appropriately, the terms in $\tilde{f}(\sign(x)\cdot |x|^q)$ will be rational in $|x_i|$ with the possible presence of $\sign(x_i)$ terms. Then, if the $\lambda_i$ terms are chosen appropriately, the rational terms will be eliminated. To account for the remaining $\sign(x)$ and $|x|$ terms, if present, an inequality of the form of Eqn.~\eqref{eq:cond4} should be imposed for each sector of the state-space -- i.e. when $x \ge 0$ use valid inequality $g(x)=x$ and when $x \le 0$ use valid inequality $g(x)=-x$. See the numerical examples for illustration.
\end{note}
%%%%%%%%%%%%%%%%%%%%%%%%%%%%%%%%%%%%%%%%%%%%%%%%%%%%%%%%%%%%%%%%%%%%%%%%%%%%%%%%%%%%%%%%%%%%%%%%%%%%%%%%%%%%%%%%%%%%%%
\section{Numerical Examples}\label{sec:numerical} \vspace{-2mm}
%This section presents numerical examples of finite-time stability analysis using SoS programming, building on the results from the preceding section.
To illustrate the application of Prop.~\ref{SoSApp}, we consider both scalar and  multivariate vector fields. The SoS conditions in both cases are enforced using SOSTOOLS~\citep{sostools}.
\begin{exmp}\label{ex:1}
Consider the scalar system given by
\begin{equation}\label{ex1}
\dot x(t)=f(x(t)) := -\sign(x(t))|x(t)|^{2/3}.
\end{equation}
We now apply Prop.~\ref{SoSApp} on domain $\mcl D=\R$ with $\Omega=\{x \in \mathbb{R} : x ^2 \leq 2 \}$ where $q=3, \lambda=1, \tau=2, d=2,  p=3,$ and $g(x)=2-x^2$. Then
\begin{align*}\label{ex2}
\tilde f(\sign(x)|x|^{3})&= \frac{1}{3}f(\sign(x)|x|^{3})|\sign(x)|x|^{3}|^{-\frac{2}{3}}\\
&=-\frac{1}{3}\sign(x) x^2x^{-2}=-\sign(x).
\end{align*}
Then condition~\eqref{eq:cond4} becomes
\[
	\nabla_x V(x)x - \mu |x|^{4}- t(x)(2-x^2)  \in \Sigma_s
\]
where we choose $t \in \Sigma_s$ to be of degree 2. Choosing $k=7.10$, we bisect on $\mu$ to find the optimal $\mu^*=11.99$. The resulting $V$ is given by $V(x) = 5.42 x^4$. Prop.~\ref{SoSApp} now implies that for $\tilde \mu=\mu/k=1.69$, $\gamma=p/2d=.75$, $\tilde{\Omega} := \{z \in \mathcal{D}: \sign(z) |z|^{1/3} \in \Omega \}=\{x\,:\, |x|\le 8\}$, $\tilde V(z) = 5.42 |z|^{4/3}$ that the settling-time function is bounded by $T(z)\le 2.37 \tilde V^{1/4}= 3.61 |z|^{1/3}$ for any $z$ such that $\{x\,:\,\tilde V(x)\le \tilde V(z)\}\subset \tilde \Omega$. Hence if $|z|^2\le 8$, $z \in \tilde \Omega$, and the sublevel set $\{x\,:\,\tilde V(x)\le \tilde V(z)\}$ is simply $\{x\,:\, |x|\le |z|\} \subset \tilde \Omega$, then the settling-time function is valid for any $|z|^2\le 8$. For example, if we take initial state $z=1.2 \in \tilde \Omega$ the corresponding settling time is $3.61 |z|^{1/3}=3.84$.
Using this initial condition, the estimated settling time obtained from numerical simulation using MATLAB's \texttt{ode-23} solver is $3.14 \, \text{s} $. %The results of this simulation can be found in Fig.~\ref{fig:est}.
\end{exmp}
\begin{comment}
\begin{figure}[htbp]
\begin{subfigure}[b]{0.48\textwidth}
	\centering
	\includegraphics[width=\textwidth]{MultVar10}
	\caption{Numerical Simulation of Eqn.~\eqref{ex1} with initial condition $x(0)=1.2$. Settling time is estimated as $t_f=3.14s$.}
	\label{fig:est}
\end{subfigure}
\hfill
\centering
\begin{subfigure}[b]{0.48\textwidth}
	\centering
	\includegraphics[width=\textwidth]{MultVar9}
	\caption{Numerical Simulation of Eqn.~\eqref{ex3} with initial condition $x_1(0)=1.3$, $x_2(0)=.8$. Settling time is estimated as $t_f=3.03s$.}
	\label{fig:est1}
\end{subfigure}
\caption{Semilog plots of the evolution of the norm of the state for simulation of System~\eqref{ex1} (a) and System~\eqref{ex3} (b) with specified initial conditions. Settling time is taken as the smallest $t$ for which $\norm{x(t)}\le 10^{-6}$.}
\label{fig:combined_est}
\end{figure}
\end{comment}

\begin{exmp}\label{ex:2}
	Consider the 2-state system given by
	\begin{equation}\label{ex3}
\bmat{\dot x_1(t)\\\dot x_2(t)}=f(x(t)):=	\bmat{		-\sign(x_1(t)) |x_1(t)|^{1/2} + x_2(t)^{1/3} \\
			- x_2(t)^{1/3} }
	\end{equation}
%	\begin{equation}\label{ex3}
%\bmat{\dot x_1(t)\\\dot x_2(t)}=	\underbrace{\bmat{		-\sign(x_1(t)) |x_1(t)|^{1/2} + \sign(x_2(t))|x_2(t)|^{1/3} \\
%			- \sign(x_2(t)) |x_2(t)|^{1/3} }}_{f(x(t))}
%	\end{equation}
  We now apply Prop.~\ref{SoSApp}  on the domain $\mathcal{D} \subset \mathbb{R}^2$, with $\Omega=\{x \in \mathcal{D} : \|x\|^2 \leq  3 \}$ where $q=[q_1, q_2]=[2, 3]$, $p=4$, $d=3$, $\tau=2$, $g_1(x)=3-(x_1^2 + x_2^2)$, and $\lambda_1=\lambda_2 =2$ so that $ \prod_{i=1}^2 |x_i|^{\lambda_i}  = x_1^2 x_2^2$ . Then
  \(
   \tilde f_1(\sign(x) \cdot |x|^{q}) = \half f_1(\sign(x) \cdot |x|^{q})(|x_1|^{2})^{-\half}
  \)
	and since $(\sign(x_2) |x_2|^3)^{1/3}=x_2$, $\sign(\sign(x_1) |x_1|^2)=\sign(x_1)$ and $ |\sign(x_1) |x_1|^2|^{\half}=|x_1|$, we have:
	\(
		\tilde f_1(\sign(x)\cdot |x|^{q}) =-\half (x_1 - x_2 ) \frac{1}{|x_1|} .
	\)
	Furthermore,
\[
		\tilde f_2(\sign(x) \cdot |x|^{q}) =-\frac{1}{3}  x_2 |\sign(x_2)|x_2|^3|^{-2/3}  = -\frac{1}{3} \frac{1}{x_2} .
\]
%Hence, \vspace{-3mm}
%	\begin{equation*}
%		\tilde f(\sign(x_1) |x_1|^{2},\sign(x_2) |x_2|^{3}) =
%		\bmat{
%			-\frac{1}{2} ( x_1 - x_2 ) \frac{1}{|x_1|}\\
%			-\frac{1}{3}\frac{1}{x_2} }.
%	\end{equation*} \vspace{-3mm}
	Since
	\begin{align*}
		&\nabla_x V(x)^T \tilde f (\sign(x)|x|^{q}) \prod_{i=1}^n |x_i|^{2\lambda_i}   \\
		&= -\frac{1}{2} \partial_{x_1} V(x)  (x_1 - x_2 ) \frac{1}{|x_1|} x_1^2 x_2^2  -\frac{1}{3} \partial_{x_2} V(x) \frac{1}{x_2} x_1^2 x_2^2 \\
		&= -\frac{1}{2} \partial_{x_1} V(x)  (x_1|x_1| x_2^2- |x_1| x_2^3) - \frac{1}{3} \partial_{x_2} V(x) x_1^2 x_2,
	\end{align*}
	condition~\eqref{eq:cond4} becomes
	\begin{multline*}
		\frac{1}{2} \partial_{x_1} V(x)  (x_1|x_1| x_2^2- |x_1|x_2^3) + \frac{1}{3} \partial_{x_2} V(x) x_1^2 x_2  \\
		- \mu \|x \|^p x_1^2x_2^2-t(x)(3-x_1^2 - x_2^2)   \in \Sigma_s.
	\end{multline*}
	where we choose $t \in \Sigma_s$ to be of degree $4$.
	Due to the presence of absolute values, condition~\eqref{eq:cond4} is applied separately on two distinct semi-algebraic sets. For $x>0$, the absolute value reduces to $|x|=x$, and the constraint is enforced by augmenting the semialgebraic set with $g_2(x)=x$. Condition~\eqref{eq:cond4} then becomes
	\begin{multline*}
		\frac{1}{2} \partial_{x_1} V(x)  (x_1^2 x_2^2- x_1x_2^3) + \frac{1}{3} \partial_{x_2} V(x) x_1^2 x_2    \\
		- \mu \|x \|^p x_1^2x_2^2- t_1(x)(3-x_1^2 - x_2^2)- v_1(x) x  \in \Sigma_s,
	\end{multline*}
	where $v_1 \in \Sigma_s$ (degree $4$).
	 For $x<0$, we likewise modify the the original semialgebraic by appending $g_2(x)=-x$. The second instance of Condition~\eqref{eq:cond4} then becomes
	 \begin{multline*}
	 	\frac{1}{2} \partial_{x_1} V(x)  (-x_1^2 x_2^2+ x_1x_2^3) + \frac{1}{3} \partial_{x_2} V(x) x_1^2 x_2    \\
	 	- \mu \|x \|^p x_1^2x_2^2- t_2(x)(3-x_1^2 - x_2^2)+ v_2(x) x \in \Sigma_s,
	 \end{multline*}
	 where $v_2 \in \Sigma_s$ (degree $4$).
%	The semialgebraic sets $g_1, g_2$ partition the domain $\mathcal{D}$ of $f$ into regions where the absolute value behaves linearly.
Choosing $k=.2$, we bisect on $\mu$ to find the optimal $\mu^*=.16$. The resulting $V$ is given by
\begin{multline*}
V(x) = 0.0160 x_1^6 - 5.1188 \times 10^{-7} x_1^5 x_2 + 0.3493 x_1^4 x_2^2 \\
- 0.1610x_1^3 x_2^3 + 0.2107 x_1^2 x_2^4 - 0.0101 x_1 x_2^5 + 0.1993 x_2^6.
\end{multline*}
Prop.~\ref{SoSApp} now implies that for $\tilde \mu=\mu/k=.802$, $\gamma=p/2d=2/3$, $\tilde{\Omega} := \{z_1, z_2 \in \mathcal{D}:(\sign(z_1) |z_1|^{1/2}, \sign(z_2) |z_2|^{1/3}) \in \Omega \}=\{(x_1,x_2):\, |x_1|+|x_2|^{2/3}\le 3 \}$,
\begin{align*}
	\tilde V(z) = &\; 0.016 z_1^3 - 5.12 \cdot 10^{-7} \, \sign(z_1) |z_1|^{5/2} z_2^{1/3} \\
	&\; + 0.35 z_1^2 z_2^{2/3} - 0.161 \, \sign(z_1) |z_1|^{3/2} z_2 \\
	&\; + 0.21 z_1 |z_2|^{4/3} - 0.01 \, \sign(z_1) |z_1|^{1/2} z_2^{5/3} + 0.2 z_2^2
\end{align*}
that the settling-time function is bounded by $T(z)\le 3.7405 \tilde V^{1/3}$ for any $z$ such that $\{x\,:\,\tilde V(x)\le \tilde V(z)\}\subset \tilde \Omega$. For example, if we take initial state $z_1(0)=1.3$ and  $z_2(0)=0.8$, then $\tilde V(z(0))=.676$ and it can be shown (using an auxiliary SoS program) that $\tilde V(x)\le .674$ implies $\norm{x}^2\le 3$. Thus, the settling-time function is valid for this initial condition and  is bounded by $3.28$.
Using this initial condition, the estimated settling time obtained from numerical simulation using MATLAB's \texttt{ode-23} solver is $3.03 \, \text{s} $. %The results of this simulation can be found in Fig.~\ref{fig:est1}.
\end{exmp}
%%%%%%%%%%%%%%%%%%%%%%%%%%%%%%%%%%%%%%%%%%%%%%%%%%%%%%%%%%%%%%%%%%%%%%%%%%%%%%%%%%%%%%%%%%%%%%%%%%%%%%%%%%%%%%%%%%%%%%%%%%%%%%%%%%%%%%%%%%%%%%%%%%%%%%%%%%%%%%%%%%%%%%%%%%%%%%%%%%%%%%%%%%%%%%%%%%%%%%%%%%%%%%%%%%%%%%%
\section{Conclusion}\vspace{-2mm}
We have proposed a Sum of Squares (SoS) framework for testing finite-time stability and bounding the associated settling-time function. For finite-time stable systems, the use of polynomial Lyapunov functions results in stability conditions with fractional exponents and signum functions -- conditions which cannot be enforced using traditional SoS programming. We resolve this issue by applying a coordinate transformation which allows polynomial inequalities and SoS programming to be used to find non-polynomial Lyapunov functions and enforce the resulting finite-time stability conditions. By solving the resulting polynomial inequalities, we obtain accurate bounds on the settling time and region of attaction. Numerical examples are used to illustrate the approach and demonstrate accuracy of the resulting bounds on settling-time. %These results have the potential to be used in sliding mode control to guarantee finite convergence to the sliding surface in a sufficiently short time.
%%%%%%%%%%%%%%%%%%%%%%%%%%%%%%%%%%%%%%%%%%%%%%%%%%%%%%%%%%%%%%%%%%%%%%%%%%%%%%%%%%%%%%%%%%%%%%%%%%%%%%%%%%%%%%%%%%%%%%%%%%%%%%%%%%%%%%%%%%%%%%%%
\bibliography{ifacconf}             % bib file to produce the bibliography

\end{document}